\magnification\magstep1
\baselineskip = 18pt
\def\vp{\varepsilon}
\def\n{\noindent}
\def\eop{\vrule height7pt width7pt depth0pt}
\overfullrule = 0pt

\def \Ball {{\rm Ball\,}}
\def \Sphere {{\rm Sphere\,}}
\def \e{\epsilon }
\def \a{\alpha }

\def \d{\delta}

\def \lam{\lambda}

\def \fn{\{f_n\}_{n=1}^\infty}
\def \hn{\{h_n\}_{n=1}^\infty}

\def \wn{\{w_n\}_{n=1}^\infty}
\def \xalphastar{\{ x_{\alpha}^* \}}
\def \yalphastar{\{ y_{\alpha}^* \}}
\def \spa{{\rm span\,}}
\def \sign{{\rm sign\,}}
\def\ss#1{_{\lower3pt\hbox{$\scriptstyle #1$}}}

\centerline{Extension of Operators from Weak$^*$-closed Subspaces of
$\ell_1$}
 \centerline{into $C(K)$ Spaces}\medskip
\centerline{W.B.~Johnson$^{*\dag}$ and M.~Zippin$^{\dag\ddag}$}\vskip1in
\centerline{\bf Abstract}

It is proved that every operator from a weak$^*$-closed subspace of
$\ell_1$ into a space $C(K)$ of continuous functions on a compact Hausdorff
space $K$ can be extended to an operator from $\ell_1$ to $C(K)$.

\vfill

\n Mathematics Subject Classification.  Primary 46E15, 46E30;
Secondary 46B15.

\bigskip\bigskip\bigskip

\n $^*$Supported in part by NSF DMS-9306376.

\n $^{\dag}$Supported in part by a grant of the U.S.-Israel Binational
Science Foundation.

\n $^{\ddag}$Participant at Workshop in Linear Analysis and Probability,
NSF DMS-9311902\eject

\n {\bf 1. Introduction.}

This work is part of an effort to characterize
those subspaces $E$ of a Banach space $X$ for which the pair $(E,X)$ has
the following:

\proclaim  Extension Property. (E.P., in short):\ Every (bounded,
linear) operator
$T$ from
$E$ into any $C(K)$ space $Y$ has an extension ${\bf T}\colon \ X\to Y$.

\n There is a quantitative version of the E.P.: \ for any $\lambda \ge 1$ we
 say that the pair $(E,P)$  has the $\lambda$-EP  if for every
$T\colon \ E\to Y$ there is an extension ${\bf T}\colon \ X\to Y$ with
$\|{\bf T}\|\le \lambda\|T\|$. It is easy to see that if $(E,X)$ has the
E.P., then it has the $\lambda$-E.P.\ for some $\lambda$.

It is known [Zip] that for each $1<p<\infty$ and every subspace $E$ of
$\ell_p$, $(E,\ell_p)$ has the 1-E.P., while for $F\subset c_0$, $(F,c_0)$
has the $(1+\vp)$-E.P.\ for  every $\vp>0$ [LP]. However, there is a
subspace $F$ of $c_0$ for which $(F,c_0)$ does not have the 1-E.P. [JZ2].
If $E$ itself is a $C(K)$ space then, clearly, $(E,X)$ has the E.P.\ if and
only if $E$ is complemented in $X$. It follows from [Ami] that $C(K)$ has a
subspace $E$ for which $(E, C(K))$ does not have the E.P.\ if $K$ is any
compact metric space whose $\omega$-th derived set is nonempty (which is
equivalent [BePe] to saying that $C(K)$ is not isomorphic to $c_0$).

Since every separable Banach space is a quotient of $\ell_1$, the following
fact demonstrates the important r\^ole of the space $\ell_1$ in extension
problems.

\proclaim Proposition 1.1. Let $E$ be a subspace of a Banach space $X$ and 
let $Q$ be an
operator from $Z$ onto $X$ so that $\|Q\|=1$ and $Q \hbox{ Ball } Z \supset
\delta \hbox{ Ball } X$. If $(Q^{-1}E, \ell_1)$ has the $\lambda$-E.P.\ then
$(E,X)$ has the $\lambda/\delta$-E.P.

\n {\bf Proof.} Let $T$ be an operator from $E$ into any $C(K)$ space $Y$.
Consider the operator $S = TQ\colon\ Q^{-1}E\to Z$. If ${\bf S}\colon \
Z\to Y$ extends $S$ then since $\bf S$ vanishes on $\hbox{ker } Q, {\bf S}$
induces an operator $\widetilde S$ from $X\sim Z/\hbox{ker } Q$ into $Y$ so
that $\widetilde S Q={\bf S}$ and $\|\widetilde S\| \le \|{\bf
S}\|/\delta$. $\hfill\eop$

An immediate consequence of Proposition~1 is that $\ell_1$ contains a
subspace $F$ for which $(F,\ell_1)$ does not have the E.P. Indeed, if $E$
denotes an uncomplemented supspace of $C[0,1]$ which is isomorphic to
$C[0,1]$ ([Ami]) and if $Q\colon \ \ell_1\to C[0,1]$ is
a quotient map and
$F = Q^{-1}E$, then $(F,\ell_1)$ does not have the E.P. The main purpose of
this paper is to prove the following.

\proclaim Theorem. Let $\{X_n\}^\infty_{n= 1}$ be finite dimensional and
let $E$ be a weak$^*$-closed subspace of $X = (\sum X_n)_1$, regarded as
the dual of $X_* = (\sum X^*_n)_{c_0}$. Then $(E,X)$ has the E.P. Moreover,
if $E$ has the approximation property, then $E$ has the $(1+\vp)$-E.P.\ for
every $\vp>0$.

We know very little about the extension problem for general pairs $(E,X)$.
However, the theorem  makes the following small contribution in the
general case.

\proclaim Corollary 1.1. Let $E$ be a subspace of the separable space $X$.
Assume that there is a weak$^*$-closed subspace $F$ of $\ell_1$ such that
$X/E$ is isomorphic to $\ell_1/F$. Then $(E,X)$ has the E.P.

\n {\bf Proof.} Let $Q\colon \ \ell_1\to X$ and $S\colon \ X\to X/E$ be
quotient maps. Theorem~2 of [LR] implies that there is an automorphism of
$\ell_1$ which maps $Q^{-1}E = \ker(SQ)$ onto $F$. Since $(F,\ell_1)$ has
the E.P.\ by the theorem, so does the pair $(Q^{-1}E,\ell_1)$. It follows
from Proposition~1 that $(E,X)$ has the E.P.$\hfill\eop$

We use standard Banach space theory notation and terminology, as may be
found in [LT1], [LT2].

\vfill\eject

\n {\bf 2. Preliminaries.}

Let $E$ be a subspace of $X$, $\lambda\ge 1$, and $0<\vp<1$. Given an
operator $S\colon \ E\to Y$ we say that the operator $T\colon \ X\to Y$ is
a $(\lambda,\vp)$-{\sl approximate extension\/} of $S$ if $\|T\| \le
\lambda \|S\|$ and
$$\|S-T|_E\| \le \vp\|S\|.$$
Our first observation is that the existence of approximate extensions
implies the existence of extensions.

\proclaim Lemma 2.1. Let $E$ be a subspace of $X$ and assume that each
operator $S\colon \ E\to Y$ has a $(\lambda,\vp)$-approximate extension.
Then the pair $(E,X)$ has the $\mu$-E.P.\ with $\mu\le
\lambda(1-\vp)^{-1}$.

\n {\bf Proof.} Put $S_1=S$ and let $T_1$ be a $(\lambda,\vp)$-approximate
extension of $S_1$. Then $\|T_1\| \le \lambda\|S_1\| = \lambda\|S\|$ and
$\|S_1-T_1|_E\|\le \vp\|S\|$. Construct by induction sequences of
operators $\{S_n\}^\infty_{n=1}$ from $E$ into $Y$ and
$\{T_n\}^\infty_{n=1}$ from $X$ into $Y$ such that for each $n\ge 1$
$S_{n+1} = S_n-T_n|_E$ and $T_{n+1}$ is a $(\lambda,\vp)$-approximate
extension of $S_{n+1}$. Then, by definition, $\|T_n\|\le \lambda\|S_n\|$
and $\|S_{n+1}\|\le \vp\|S_n\|$ for every $n\ge 1$. It follows that
$\left\|S -\sum\limits^n_{i=1} T_i|_E\right\|\le \vp^n\|S\|$ and $\|T_n\|
\le \lambda \vp^{n-1}\|S\|$ for all $n\ge 1$. Hence the operator $T =
\sum\limits^\infty_{i=1} T_i$ extends $S$ and $\|T\| \le
\lambda(1-\vp)^{-1}\|S\|$.$\hfill \eop$

Given a finite dimensional decomposition (FDD, in short)
$\{Z_n\}^\infty_{n=1}$ of a space $Z$, we will be interested in subspaces
 of $Z$ with FDD's which are particularly well-positioned with respect to
$\{Z_n\}^\infty_{n=1}$.

\n {\bf Definition.} Let $F\subset Z$ and let $\langle
F_n\rangle^\infty_{n=1}$
be an FDD for $F$. We say that $\{F_n\}^\infty_{n=1}$ is {\sl alternately
disjointly supported\/} with respect to $\{Z_n\}^\infty_{n=1}$ if there
exist integers $1=k(1) < k(2)<\cdots$ such that for each $n\ge 1$, $F_n
\subset Z_{k(n)} + Z_{k(n)+1} +\cdots+ Z_{k(n+2)-1}$. 

An important property
of an alternatively disjointly supported FDD is that if
$\{n(j)\}^\infty_{j=1}$ is any increasing sequence of integers and if we
drop $\{F_{n(j)}\}^\infty_{j=1}$,  then the remaining $F_n$'s can be grouped
into blocks $\widetilde F_j = \sum\limits^{n(j+1)-1}_{i=n(j)+1}F_i$ which
form
an FDD that is disjointly supported on the $\{Z_n\}^\infty_{n=1}$; more
precisely, with the above notation,
$$\widetilde F_j \subset \sum^{k(n(j+1)+1)-1}_{m=k(n(j)+1)} Z_m\qquad
\hbox{for all}\quad j\ge 1.$$
We will show that for certain subspaces of a dual space with an FDD, a
given FDD can be replaced by one which is alternately disjointly supported.

We first need the following main tool:

\proclaim  Proposition 2.1. Let $\{X_n\}^\infty_{n=1}$ be a shrinking FDD
for $X$, let $Q$ be a quotient mapping of $X$ onto $Y$ and suppose that
$\{\widetilde E_n\}^\infty_{n=1}$ is an FDD for $Y$. Then there are a
blocking $\{E'_n\}^\infty_{m=1}$ of $\{\widetilde E_n\}^\infty_{n=1}$, an
FDD $\{W_n\}^\infty_{n=1}$ of $X$ which is equivalent to
$\{X_n\}^\infty_{n=1}$, and $1=k(1) < k(2)<\cdots$ so that for each $n$ and
each $k(n) \le j < k(n+1)$, $QW_j \subset E'_n + E'_{n+1}$. Moreover, given
$\vp>0$, $\{E'_n\}^\infty_{n=1}$ and $\{W_n\}^\infty_{m=1}$ can be chosen
so that there is an automorphism $T$ on $X$ with $\|I-T\| < \vp$ and $TX_n
=W_n$ for all $n$.

\n {\bf Proof.} In order to avoid complicated notation we shall prove the
statement for the case where, for every $n\ge 1$, $X_n$ (and hence also
$W_n$) is one dimensional.
The same arguments, with only obvious modifications yield the FDD case.
(Actually, in the proof of the theorem, only the basis case of
Proposition~2 is needed. Indeed, in Step~3 of the proof of the theorem, one
can replace $E$ by $E_1 \equiv	E\oplus_1 (\sum G_n)_1$ and $X$ by
$X_1 = X\oplus_1(\sum G_n)_1$, where $\{G_n\}^\infty_{n=1}$ is a sequence
which is dense in the sense of the Banach-Mazur distance in the set of all
finite dimensional spaces, and use the fact [JRZ], [Pel] that $E_1$ has a
basis. In fact, this trick is used in a different way for the proof of the
``moreover'' statement in the theorem.)

So assume that $X$ has a normalized shrinking basis $\{x_n\}^\infty_{n=1}$
with  biorthogonal functionals $\{f_n\}^\infty_{n=1}$; we
are looking for an equivalent basis $\{w_n\}^\infty_{n=1}$ of $X$ for which
the statement holds. First we perturb the basis for $X$ to get another basis
whose images under $Q$ are supported on finitely many of the $\widetilde
E_n$'s. This step does not require the hypothesis that
$\{x_n\}^\infty_{n=1}$ be shrinking.

For each $n\ge 1$ let $\widetilde Q_n$ be the FDD's natural projection from
$Y$ onto $\widetilde E_1 + \widetilde E_2 +\cdots+ \widetilde E_n$. Let
$1>\epsilon>0$ and set $C=\sup_n\|f_n\|$. Choose $p_1 < p_2<\cdots$
so that for each $n$, $\|Qx_n-\widetilde Q_{p_n}Qx_n\| < \epsilon
C^{-1}2^{-n}$. Since $Q$ is a quotient mapping, there is for each $n$ a
vector $z_n$ in $X$ with $\|z_n\| < \epsilon C^{-1}2^{-n}$ and $Qz_n = Qx_n
- \widetilde Q_{p_n}Qx_n$. Let $y_n = x_n-z_n$, so that $Qy_n$ is in
$\widetilde E_1 +\cdots \widetilde Ep_n$. It is standard to check that
$\{y_n\}^\infty_{n=1}$ is equivalent to $\{x_n\}^\infty_{n=1}$. Indeed,
define an operator $S$ on $X$ by $Sx = \sum\limits^\infty_{n=1} f_n(x)z_n$.
Then $\|S\|<\epsilon$ and $Sx_n = z_n$, so $I-S$ is an isomorphism from $X$
onto $X$ which maps $x_n$ to $y_n$.

Define a blocking $\{E_n\}^\infty_{n=1}$ of $\{\widetilde
E_n\}^\infty_{n=1}$ by $E_n = \widetilde E_{p_{n-1}+1} + \cdots +
\widetilde E_{p_n}$ (where $p_0\equiv 0$). Then for each $n$, $Qy_n$ is in
$E_1+\cdots + E_n$.

Let $Q_n$ be the basis projection from $Y$ onto $E_1 +\cdots+ E_n$, $P_n$
the basis projection from $X$ onto $\hbox{span}\{y_1,\ldots, y_n\}$, and
set $C_1 = \sup_n\|P_n\|$. Since $\{y_n\}^\infty_{n=1}$ is shrinking,
$\lim\limits_{m\to \infty}\|Q_nQ(I-P_m)\|=0$. Since $Q$ is a quotient
mapping, for each $n$ there exists a mapping $T_n$ from $E_1+\cdots+ E_n$
into $X$ so that $QT_n$ is the identity on $E_1+\cdots+ E_n$. Set $M_n =
\|T_n\|$, let $1>\epsilon>0$, and recursively choose $0 = k(0) < 1 = k(1)
<k(2)<\cdots$ so that for each $n$, $\|Q_{k(n)}Q(I-P_{k(n+1)-1})\| <
(2C_1 M_{k(n)}))^{-1}2^{-n}\epsilon$. Setting $w_j  = y_j
-T_{k(n)}Q_{k(n)}Qy_j$ for $k(n+1) \le j < k(n+2)$, we see that $Qw_j$ is
in $E_{k(n)+1}+\cdots+ E_{k(n+2)}$ when $k(n+1)\le j<k(n+2)$.

The desired blocking of $\{\widetilde E_n\}^\infty_{n=1}$ is defined by
$E'_n = E_{k(n-1)+1} + E_{k(n-1)+2} +\cdots+ E_{k(n)}$, but it remains to
be seen that $\{w_n\}^\infty_{n=1}$ is a suitably small perturbation of
$\{y_n\}^\infty_{n=1}$. The inequality $\|Q_{k(n)}Q(I-P_{k(n+1)-1})\| <
(2C_1 M_{k(n)})^{-1}2^{-n}\epsilon$ implies, by composing on the right with
$P_{k(n+2)-1}$, that $\|Q_{k(n)}Q(P_{k(n_2)-1}-P({k(n+1)-1})\| <
(2M_{k(n)})^{-1}2^{-n}\epsilon$. Thus if we define an operator $V$ on $X$
by $Vx = \sum^\infty_{n=0} T_{k(n)} Q_{k(n)}Q(P_{k(n+2)-1} -
P_{k(n+1)-1})x$, we see that $\|V\| < \epsilon$ and hence $T\equiv I-V$ is
invertible. But for $k(n+1) \le j < k(n+2)$, $Vy_j = T_{k(n)}Q_{k(n)}Qx_j$;
that is, $Ty_j = w_j$. \hfill\eop

Using a duality argument we get from Proposition~2.1 the following.

\proclaim Corollary 2.1. Let $\{Z_n\}^\infty_{n=1}$ be an $\ell_1$-FDD for
a space $Z$. Regard $Z$ as the dual of the space $Z_* = (\sum Z^*_n)_{c_0}$
and let $F$ be a weak$^*$-closed subspace of $Z$ with an FDD. Then $Z$ and
$F$ have $\ell_1$-FDD's $\{V_n\}^\infty_{n=1}$ and $\{U_n\}^\infty_{n=1}$,
respectively, so that $\{U_n\}^\infty_{n=1}$ is alternately disjointly
supported with respect to $\{V_n\}^\infty_{n=1}$. Moreover, given $\vp>0$,
$\{V_n\}^\infty_{n=1}$ can be chosen so that for some blocking
$\{Z'_n\}^\infty_{n=1}$ of $\{Z_n\}^\infty_{n=1}$, there is an automorphism
$T$ of $Z_*$ with $\|I-T\|<\vp$ and $TZ'_n = V_n$ for all $n\ge 1$.

\n {\bf Proof.} Being weak$^*$-closed, $F$ has a predual $F_* = Z_*/F_\bot$
which is a quotient space of $Z_*$. By [JRZ], $F_*$ has a
shrinking FDD and consequently, by Theorem~1 of [JZ2], $F_*$ has a
shrinking $c_0$-FDD $\{\widetilde E_n\}^\infty_{n=1}$. Let $Q\colon \
Z_*\to F_*$ be the quotient mapping.  By Proposition~2.1 there are a
blocking $\{E'_n\}^\infty_{n=1}$ of $\{\widetilde E_n\}^\infty_{n=1}$, an
FDD $\{W_n\}^\infty_{n=1}$ of $Z_*$ which is equivalent to
$\{Z^*_n\}^\infty_{n=1}$, even the image of $\{Z^*_n\}^\infty_{n=1}$ under
some automorphism on $Z_*$ which is arbitrarily close to $I_{Z_*}$, and
$1=k(1) < k(2) <\cdots$ so that for each $n$ and $k(n) \le j<k(n+1)$, $QW_j
\subset E'_n +E'_{n+1}$. The equivalence implies that
$\{W_n\}^\infty_{n=1}$ is a $c_0$-FDD and, being a blocking of a $c_0$-FDD,
$\{E'_n\}^\infty_{n=1}$ is a $c_0$-FDD. Let $\{V_n\}^\infty_{n=1}$ (resp.\
$\{U_n\}^\infty_{n=1}$) be the dual FDD of $\{W_n\}^\infty_{n=1}$ (resp.
$\{E'_n\}^\infty_{n=1}$) for $Z$ (resp.\ $F$). Then $\{V_n\}^\infty_{n=1}$
is an $\ell_1$-FDD for $Z$ and $\{U_n\}^\infty_{n=1}$ is an $\ell_1$-FDD for
$F$. Moreover, suppose that $u$ is in $U_n$ and $w_j$ is in $W_j$, where
either $j<k(n)$ or $j\ge k(n+2)$. Let $m$ be the integer for which $k(m) \le
j<k(m+1)$. Then either $m<n$ or $m>n+1$ hence $n\ne m$ and $n\ne m+1$.
Then \ $Qw_j \in E'_m + E'_{m+1}$, hence $u(w_j) = \langle u,Qw_j\rangle
= 0$. This proves that $U_n$ is supported on \ \  
$\sum\limits^{k(n+2)-1}_{j=k(n)}V_j$.$\hfill \eop$
\vfill\eject

\n {\bf 3. Proof of the Theorem.}

The proof consists of four parts, the first three of which are essentially
simple special cases of the theorem.

\proclaim Step 1. $E$ has an  FDD \ $\{E_n\}^\infty_{n=1}$ \ with \ \ 
$ E_n\subset X_n$ \  for all \ $n$.

\n {\bf Proof.} Let $Y	= C(K)$ and let $S\colon \ E\to Y$ be any operator.
Using the ${\cal L}_{\infty,1+\vp}$-property of $Y$ (or see Theorem~6.1 of
[Lin]), one sees that the finite rank operator $S|_{E_n}$ has an
extension $S_n\colon \ X_n\to Y$ with
$\|S_n\|\le (1+\vp)\|S_n\|$. Define the extension $\bf S$ of $S$ by ${\bf
S}\left(\sum\limits^\infty_1 x_n\right) = \sum\limits^\infty_{n=1}
S_nx_n$. Since \ $\{X_n\}^\infty_{n=1}$ \ is an exact
$\ell_1$-decompostion, it follows that $\|{\bf S}\| \le
(1+\vp)\|S\|$.

\proclaim Step 2. $E$ has an $\ell_1$-FDD $\{E_n\}^\infty_{n=1}$ which is
alternately disjointly supported with respect to $\{X_n\}^\infty_{n=1}$.

\n {\bf Proof.} Given $\delta>0$, let $1<(1+\vp)(1-\vp)^{-1}<1+\delta$ and
choose an integer $N>(1+\vp)M\vp^{-1}$ where $M$ is the constant of the
$\ell_1$-FDD $\{E_n\}^\infty_{n=1}$; that is, the constant of equivalence
of \ $\{E_n\}^\infty_{n=1}$ \ to the natural $\ell_1$-FDD for $(\sum
E_n)_1$. Let
$Y =C(K)$ and let
$S\colon
\ E\to Y$ be an operator  with $\|S\|=1$. For each $1\le j \le N$ let 
$$Z_j =
\overline{\rm span}\{E_i \, \colon \ i\ne kN+j, \ k=0,1,2, \ \ldots\}.$$
 Each
subspace $Z_j$ has a natural $\ell_1$-FDD  which is disjointly supported
with respect to $\{X_n\}^\infty_{n=1}$ because $\{E_n\}^\infty_{n=1}$
is alternately disjointly supported with respect to $\{X_n\}^\infty_{n=1}$.
By Step~1, $S|_{Z_j}$ has an extension $T_j\colon \ X\to Y$ with 
$$
\|T_j\|
\le (1+\vp)\|S_j\| \le (1+\vp) \|S\| = 1+\vp.
$$
 Define $T\colon \ Z\to Y$ by
$T = N^{-1} \sum\limits^N_{j=1} T_j$. Then $\|T\| \le (1+\vp) \|S\| =
1+\vp$. Moreover, if $e\in E_i$ and $i=kN+h$ for some $1\le h\le N$, then
$T_je = S_je=Se$ for all $j\ne h$ hence $T$ is ``almost'' an extension of
$S$. Indeed, \ $\|Te-Se\| = {1\over N}\|T_he-Se\| \le {2+\vp\over N}\|e\|$
\ whenever $e\in E_i$ for some $i$. Recalling that the $\ell_1$-FDD
$\{E_n\}^\infty_{n=1}$ has constant
$M$, we have that
$$\|T|_E - S\| \le M \sup_n \|T|_{E_n} - S|_{E_n}\| \le {M(2+\vp)\over N} <
\vp.$$
This proves that $T$ is an $(1+\vp,\vp)$-approximate extension of $S$ and
therefore, by Lemma~2.1, $(E,Z)$ has the $(1+\vp)(1-\vp)^{-1}$-E.P. 

\proclaim Step 3. $E$ has an FDD.

\n {\bf Proof.} By Corollary~2.1,  $X$ and $E$ have $\ell_1$-FDD's
$\{Z_n\}^\infty_{n=1}$ and $\{E_n\}^\infty_{n=1}$, respectively, with
$\{E_n\}^\infty_{n=1}$ is alternately disjointly supported with respect
to $\{Z_n\}^\infty_{n=1}$, and, by Remark~2.1, $\{Z_n\}^\infty_{n=1}$ has
constant of equivalence to $(\sum Z_n)_1$ arbitrarily close to one. Hence,
by Step~2, $(E,X)$ has the $(1+\delta)$-E.P.\ for every $\delta>0$.

This gives the ``moreover'' statement when $E$ has an FDD. When $E$ just
has the approximation property, we enlarge $X$ to $X_1 \equiv X\oplus _1
C_1$, where $C_1 = (\sum G_n)_1$ and $\{G_n\}^\infty_{n=1}$ is a sequence
of finite dimensional spaces which is dense (in the sense of the
Banach-Mazur distance) in the set of all finite dimensional spaces; and we
enlarge $E$ to $E_1 \equiv E \oplus_1 C_1$. $X_1$ is again an exact
$\ell_1$-sum of finite dimensional spaces and $E_1$ is weak$^*$-closed in
$X_1$. Moreover, since $E$ is a dual space which has the approximation
property, $E$ has the metric approximation property [LT1], and hence by
[Joh],
$E_1$ is a $\pi$-space, whence, since $E_1$ is a dual space, $E_1$ has an 
FDD by [JRZ]. Thus by Step~3,
$(E_1,X_1)$ has the $(1+\delta)$-E.P.\ for each $\delta>0$, and, therefore,
so does $(E,X)$.

\proclaim Step 4. The general case.

We start with a lemma.

\proclaim Lemma 3.1. Let $Z$ be a Banach space and let $E$ be a subspace of
$Z$. Suppose that $E$ has a subspace $F$ such that $(F,Z)$ has the
$\lambda$-E.P.\ and $(E/F, Z/F)$ has the $\mu$-E.P. Then $(E,Z)$ has the
$(\lambda+\mu(1+\lambda))$-E.P.

\n {\bf Proof.} Let $Y=C(K)$ and let $S\colon \ E\to Y$ be any operator.
Let $S_1\colon \ Z\to Y$ be an extension of $S_{|F}$ with $\|S_1\|\le
\lambda\|S\|$. The operator $W=S-S_1|_E$ from $E$ into $Y$ vanishes on $F$
and so induces an operator $\widetilde W\colon \ E/F\to Y$ in the usual
way, and $\|\widetilde W\| = \|W\| \le \|S\| + \|S_1\| \le
(1+\lambda)\|S\|$. By our assumptions, $\widetilde W$ extends to an
operator $W_1\colon \ Z/F\to Y$ with $\|W_1\| \le \mu\|\widetilde W\| \le
\mu(1+\lambda)\|S\|$. Let $Q\colon \ Z\to Z/F$ denote the quotient map.
Then $T = S_1 + W_1Q$ is the desired extension of $S$. Indeed, for every
$e \in E$
$$Te = S_1e + W_1Qe = S_1e + We  = S_1e+(S-S_1)e = Se$$
and $\|T\| \le \|S_1\| + \|W_1\| \le (\lambda+\mu(1+\lambda))\|S\|$.
$\hfill \eop$

Let us now return to the proof of the general case. Being a  weak$^*$-closed
subspace of $\ell_1$, $E$ is the dual of the quotient space \ $E_* = (\sum
X^*_n)_{c_0}/E_\top$. Our main tool in this part of the proof is
Theorem~IV.4 of [JR] and its proof. This theorem states that $E_*$ has a
subspace $V$ so that both $V$ and $E_*/V$ have shrinking FDD's. Under these
circumstances, Theorem~1 of [JZ1] implies that both $V$ and $E_*/V$ have
$c_0$-FDD's. In order to prove the theorem it suffices, in view of
Lemma~3.1, to show that both pairs $(V^\bot, X)$ and $(E/V^\bot, X/V^\bot)$
have the E.P. Now $(V^\bot,X)$ has the \ $(1+\delta)$-E.P.\ for all
$\delta>0$ by Step~3, so it remains to discuss the pair $(E/V^\bot,
X/V^\bot)$. This discussion requires some preparation and some minor
modification in the proof of Theorem~IV.4. of [JR]. We first need a known
perturbation lemma:

\proclaim Lemma 3.2. Suppose $E,F$ are subspaces of $X^*$ with $F$ norm
dense in $X^*$ and $X^*$ is separable. Then for each $\vp>0$ there is an
automorphism $T$ on $X$ so that $\|I-T\|<\vp$ and $T^*E\cap F$ is norm
dense in $T^*E$.

\n {\bf Proof.} Let $(x_n,x^*_n)$ be a biorthogonal sequence in $X\times E$
with
$\overline{\rm span}\ x^*_n = E$ (see, e.g., [Mac]) and take $y^*_n \in F$
so that $\sum\|x^*_n-y^*_n\|\, \|x_n\|<\vp$. Define $T\colon \ X\to X$ by
$$Tx = x - \sum^\infty_{n=1} \langle x^*_n-y^*_n,x\rangle
x_n.\qquad\qquad
\eop$$ 
Returning to the proof of the theorem, we may assume, in view of
Lemma~3.2, that $E\cap \hbox{span } \bigcup\limits^\infty_{n=1} X_n$ is
norm dense in
$E$. The standard back-and-forth technique [Mac] for producing biorthogonal
sequences yields a biorthogonal sequence $\{(x_n,x^*_n)\}^\infty_{n=1}
\subset X_*\times E$ with $\hbox{span}\{Qx_n\}^\infty_{n=1} = \hbox{span }
\bigcup\limits^\infty_{n=1} QX^*_n$, $\hbox{span}\{x^*_n\}^\infty_{n=1} =
E\cap \hbox{span } \bigcup\limits^\infty_{n=1} X_n$, and where $Q$ is the
quotient mapping from the predual $X_* = (\sum X^*_n)_{c_0}$ of $X$ onto
the predual $E_*$ of $E$. This means that for any $N$, $x^*_j$ is in
$\hbox{span } \bigcup\limits^\infty_{n=N} X_n$ if $j$ is sufficiently
large.

We now refer to the construction in Theorem~IV.4 of [JR] and the finite
sets $\Delta_1\subset \Delta_2\subset\cdots$ of natural numbers defined
there. From that construction, it is clear that, having defined $\Delta_n$,
the smallest element, $k(n)$, in $\Delta_{n+1}\setminus\Delta_n$ can be
as large as we desire. In particular, if $\{x^*_j\}^{\rm max\
\Delta_n}_{j=1}$ is a subset of $\hbox{span } \bigcup\limits^{m(n)}_{i=1}
X_i$, then we choose $k(n)$ large enough so that for $j\ge k(n)$, $x^*_j$
is in
$\hbox{span } \bigcup\limits^\infty_{i=m(n)+1} X_i$. Thus setting 
$$Z_n =
\hbox{span } \{ x^*_j \colon \ \  j\in \Delta_n\setminus
\Delta_{n+1}\}$$ (where $\Delta_0 \equiv \emptyset$), we have that
$\{Z_n\}^\infty_{n=1}$ is disjointly supported relative to
$\{X_n\}^\infty_{n=1}$. (In the notation above and setting $m(0) = 0$, we
have for each $n$ that
$$Z_n \subset \hbox{span}\{X_j\}^{m(n)}_{j=m(n-1)+1}.\leqno (*)$$

The subspace $V$ of $E_*$ is defined to be the annihilator of
$\left\{x^*_j\colon \ j \in \bigcup\limits^\infty_{n=1} \Delta _n\right\}$
and, as mentioned earlier, it follows from [JR] and [JZ1] that $V$ has a
$c_0$-FDD and thus $V^* = E/V^\bot$ has an $\ell_1$-FDD. It is also proved
in [JR], but is obvious from the ``extra'' we have added here, that
$\overline{\rm span}\{Z_j\}^\infty_{j=1}$ is weak$^*$-closed and hence
equals $V^\bot$. It is also obvious from $(*)$ that $X/V^\bot$ has an
$\ell_1$-FDD. Therefore, by Step~3 $(E_*/V^\bot, X/V^\bot)$ has the E.P.
$\hfill\eop$
\medskip

\n {\bf Remark.} Under the hypotheses of the theorem, we do not know
whether $(E,X)$ has the $(1+\vp)$-E.P.\ for every $\vp>0$ when $E$ fails the
approximation property. The proof we gave yields only that $(E,X)$ has the
$(3+\vp)$-E.P.\ for all $\vp>0$.

\vfill\eject

\n {\bf 4. Concluding Remarks and Problems.}

Very little is known about the Extension Property, so there is no shortage
of problems.

\proclaim Problem 4.1. If $E$ is a subspace of $X$ and $X$ is reflexive,
does $(E,X)$ have the E.P.?  What if $X$ is superreflexive?  What if $X$
is $L_p$, $1<p\not=2<\infty$?

\proclaim Problem 4.2.  If $E$ is a reflexive subspace of $X$, does
$(E,X)$ have the E.P.?  What if $E$ is just isomorphic to a conjugate
space? In the latter case, what if, in addition, $X$ is $\ell_1$?

If $E$ is a subspace of $c_0$, then $(E,c_0)$ has the $(1+\vp)$-E.P. for
every $\vp>0$ [LP] but need not have the $1$-E.P. [JZ2].  We do not know
if this phenomenon can occur in the setting of reflexive spaces:

\proclaim Problem 4.3.  If $X$ is reflexive and $(E,X)$ has the
$(1+\vp)$-E.P. for every $\vp>0$, does $(E,X)$ have the $1$-E.P.?

The following observation gives an affirmative answer to Problem 4.3 in a
special case.

\proclaim Proposition 4.1. If $X$ is uniformly smooth and $(E,X)$ has the
$(1+\vp)$-E.P. for every $\vp>0$, then  $(E,X)$ has the $1$-E.P.

\n {\bf Proof.} In preparation for the proof, we recall Proposition $2$ of
[Zip], which says:

{\sl $(E,X)$ has the $\lambda$-E.P. if and only if there exists a 
weak$^*$-continuous extension mapping from $\Ball E^*$ to $\lam \Ball X^*$; that
is, a continuous mapping \ $\phi : (\Ball E^*,{\rm weak}^*) \to (\lam\Ball
X^*,{\rm weak}^*)$ for which $(\phi e^*)|_{E}=e^*$ for every $e^*$ in
$\Ball E^*$.}

Since $X$ is uniformly smooth, given $\vp>0$ there exists $\d>0$ so that if
$x^*$, $y^*$  in $X^*$ and $x$ in $X$ satisfy \ 
$\|x^*\|=\|x\|=1=\langle x^*,x\rangle = \langle y^*,x\rangle$ \ with 
\ $\|y^*\| < 1+\d$, \ then \ $\|x^*-y^*\| < \vp$.
Letting \ $\phi_n : \Ball E^* \to (1+{1\over n})\Ball X^*$ \ be a weakly
continuous extension mapping and letting \ $f : \Sphere E^* \to \Sphere X^*$ \
be the (uniquely defined, by smoothness) Hahn-Banach extension mapping, we
conclude that
$$
\lim_{n\to\infty} \sup \{\|\phi_n(x^*)-f(x^*)\| : x^*\in\Sphere E^* \} = 0.
$$
That is, \  $\{{\phi_n|}_{\Sphere E^*}\}_{n=1}^\infty$ \ is uniformly
convergent to 
$f|_{\Sphere E^*}$.  Since each $\phi_n$ is weakly continuous, so is
$f|_{\Sphere E^*}$.

If $E$ is finite dimensional, then clearly the positively homogeneous extension
of $f$ to a mapping from $\Ball E^*$ into $\Ball X^*$ is a weakly continuous
extension mapping.  So assume that $E$ has infinite dimension.  But then 
$\Sphere E^*$ is weakly dense in $\Ball E^*$, so by the weak continuity of the
$\phi_n$'s and the weak lower semicontinuity of the norm, we have
$$
\sup \{\|\phi_n(x^*)-\phi_m(x^*)\| : x^*\in \Ball E^* \} = 
\sup \{\|\phi_n(x^*)-\phi_m(x^*)\| : x^*\in \Sphere E^* \},
$$
which we saw tends to zero as $n$, $m$ tend to infinity.  That is,
$\{\phi_n\}_{n=1}^\infty$ is a uniformly Cauchy sequence of weakly continuous
functions and hence its limit is also weakly continuous.\hfill\eop

It is apparent from the proof of Proposition 4.1 that the $1$-E.P. is fairly
easy to study in a smooth reflexive space $X$ because every extension
mapping from $\Ball E^*$ to $\Ball X^*$ is, on the unit sphere of $E^*$,
the unique Hahn-Banach extension mapping.  Let us examine this situation a
bit more in the general case.  Suppose $E$ is a subspace of $X$ and let
$A(E)$ be the collection of all norm one functionals in $E^*$ which attain
their norm at a point of
$\Ball E^*$. The Bishop-Phelps theorem [BP], [Die]  says that $A(E)$ is
norm dense in $\Sphere E^*$, hence, if $E$ has infinite dimension, 
$A(E)$ is weak$^*$ dense in $\Ball E^*$.  Therefore $(E,X)$ has the $1$-E.P. if
and only if there is a weak$^*$ continuous Hahn-Banach selection mapping $\phi :
A(E)\to\Ball X^*$ which has a weak$^*$ continuous extension to a mapping ${\bf
\phi}$ from $\overline {A(E)}^{w^*} = \Ball E^*$ to $\Ball X^*$, since
clearly 
${\bf \phi}$ will then be an extension mapping.  The existence of ${\bf
\phi}$ is equivalent to saying that whenever $\xalphastar$ is a net in $A(E)$
which weak$^*$ converges in $E^*$, then $\{\phi x_{\a}^*\}$ weak$^*$ converges in
$X^*$ (see, for example, [Bou I.8.5]).  Now when $X$ is smooth, there is only one
mapping $\phi$ to consider, and in this case the above discussion yields the
next proposition when $\dim E=\infty$ (when $\dim E<\infty$ one extends from
$\Sphere E^*=\overline {A(E)}^{w^*}$ to $\Ball E^*$ by homogeneity).

\proclaim Proposition 4.2.  Let $E$ be a subspace of the smooth space $X$.  The
pair $(E,X)$ fails the $1$-E.P. if and only if there are nets $\xalphastar$,
$\yalphastar$ of functionals in $\Sphere X^*$ which attain their norm at points of
$\Sphere E$ and which weak$^*$ converge to distinct points $x^*$ and $y^*$,
respectively, which satisfy $x^*|_{E}=y^*|_{E}$.

An immediate, but surprising to us, corollary to Proposition 4.2 is:

\proclaim Corollary 4.1.  Let $E$ be a subspace of the smooth space $X$.  If the
pair $(E,X)$ fails the $1$-E.P., then there is a subspace $F$ of $X$ of
codimension one which contains $E$ so that $(F,X)$ fails the $1$-E.P.

\n {\bf Proof.} Get  $x^*$,  $y^*$ from
Proposition 4.2 and set   $F=\spa E\cup (\ker x^* \cap \ker y^* )$.
\hfill\eop

\proclaim Problem 4.4. Is Corollary 4.1 true for a general space $X$?

\proclaim Corollary 4.2.  For $1<p\not=2<\infty$, $L_p$ has a subspace $E$ for
which $(E,L_p)$ fails the $1$-E.P.

\n {\bf Proof.} We regard $L_p$ as $L_p(0,2)$ and make the identifications
$L_p^*=L_q=L_q(0,2)$, where $q={p\over{p-1}}$ is the conjugate index to $p$.
Let 
$$f={\bf 1}_{(0,{1\over 2})} - {\bf 1}_{({1\over 2},1)}, \qquad 
g=-2\cdot {\bf 1}_{({1\over 2},1)}-{\bf 1}_{(1,2)},
$$
regarded as elements of $L_q$, and define
$$
E=(f-g)^\bot = \{ x\in L_p(0,2) : \int_0^2 x = 0\}.
$$
Notice that $|f|^{q-1} \sign f $ is in $E$, which implies that
$1=\|f\|_q=\|f\|\ss{L_p^*}=\|f_{|E}\|\ss{E^*}$.  So $f$ and $g$ induce the same
linear functional on $E$ (we write $f|_{E}=g|_{E}$), and $f$ is the unique
Hahn-Banach extension of this functional to a functional in $L_p^*=L_q$.

\proclaim Claim.  There exists $h$ in $L_q$ supported on $[0,{1\over 2}]$ so
that \ $\int_0^2 h = 0 = \int_0^2 |g+h|^{q-1}\sign (g+h)$.

Assume the claim.  Set $\lam=\|g+h\|_q$ and let $\hn$ be a sequence of functions
which have the same distribution as $h$, are supported on $[0,{1\over 2}]$, and
are probabilistically independent as random variables o $[0,{1\over 2}]$ with
normalized Lebesgue measure.  Then $g_n\equiv \lam^{-1} (g+h_n)$ defines a
sequence on the unit sphere of $L_q(0,2)$ which converges weakly to
$\lam^{-1}g$.  Moreover, 
$|g_n|^{q-1} \sign g_n$ is in $E$, which means that as a linear functional on
$L_p$, $g_n$ attains its norm at a point on the unit sphere of $E$. In view of
Proposition 4.2, to complete the proof it suffices to find a sequence
$\fn$ on the unit sphere of $L_q$ which converges weakly in $L_q$ to
$\lam^{-1}f$ so that 
$|f_n|^{q-1} \sign f_n$ is in $E$.  This is easy:  take $w$ supported on $[1,2]$
so that
$$
\int_0^2 w =0=\int_0^2 |w|^{q-1} \sign w \ \left( = \int_0^2 |f+w|^{q-1}
\sign(f+w)\right)
$$
and $\|f+w\|^q_q = 1 = 1+ \|w\|^q_q = \lam^q$ \ (so $w$ can be a
multiple of
\ ${\bf 1}_{(1,{3\over 2})}-{\bf 1}_{({3\over 2},2)}$).  Let $\wn$ be a sequence of
functions which have the same distribution as $w$, are supported on $[1,2]$, and
are probabilistically independent as random variables on $[1,2]$.  Now set 
$f_n=\lam^{-1}(f+w_n)$.

We turn to the proof of the claim.  Fix any $0<\vp<{1\over 4}$. For
appropriate
$d$, the choice
$$
h=d(4\vp {\bf 1}_{(0,{1\over 4})} - {\bf 1}_{({1\over 2}-\vp, {1\over 2})})
$$
works.  Indeed, $\int_0^2 h=0$ no matter what $d$ is, and $gh=0$, 
so we need choose $d$ to satisfy 
$$
-\int_0^2|g|^{q-1}\sign g = \int_0^2 |h|^{q-1} \sign h.
\leqno (*)$$
The left side of $(*)$ is $2^{q-1}+1>0$, while the right side is
$|d|^{q-1} \sign g \vp^{q-1}[({1\over 4})^{2-q}-\vp^{2-q}]$, so such a choice of
$d$ is possible for $p\not= 2$.
\hfill\eop

\proclaim Problem 4.5.  If $E$ is a weak$^*$-closed subspace of
$\ell_1$, does $(E,\ell_1)$  have the $1+\epsilon$-E.P. for every
$\e>0$?

A negative answer to Problem 4.5 would be particularly interesting,
because it would justify the weird approach we used to prove the
Theorem.  However, we do not even know a counterexample to:

\proclaim Problem 4.6. If $E$ is a weak$^*$-closed subspace of
$\ell_1$, does $(E,\ell_1)$  have the $1$-E.P.?

The answer to Problem 4.6 is known to be yes for finite dimensional
$E$, [Sam1], [Sam2].

\vfill\eject

\centerline{\bf References}

\item{[Ami]} D. Amir,  {\sl Continuous function spaces with the separable
projection property,\/} {\bf Bull. Res. Council Israel 10F} (1962),
163--164.

\item{[BePe]} C. Bessaga and A. Pe\l czy\'nski, {\sl Spaces of continuous
functions IV,\/} {\bf Studia Math. 19} (1960), 53--62

\item{[BP]} E. Bishop and R. R. Phelps, {\sl A proof that every Banach
space is subreflexive,\/} {\bf Bull. AMS 67} (1961), 97--98.

\item{[Bou]} N. Bourbaki, {\sl General Topology, Part 1,\/} Addison-Wesley
(1966).

\item{[Die]} J. Diestel, {\sl Geometry of Banach spaces-selcted
topics,\/} {\bf Lecture Notes in Math. 485}   Springer-Verlag (1975).

\item{[Joh]} W. B. Johnson, {\sl Factoring compact
operators,\/} {\bf Israel J. Math. 9} (1971), 337--345.

\item{[JR]} W. B. Johnson and H. P. Rosenthal, {\sl On $w^*$-basic
sequences and their applications to the study of Banach spaces,\/} {\bf
Studia Math. 43} (1972), 77--92.

\item{[JRZ]} W. B. Johnson, H. P. Rosenthal, and M. Zippin, {\sl On bases,  
finite dimensional decompositions, and weaker structures in Banach  
spaces,}\/ {\bf  Israel J. Math. 9} (1971), 488--506.

\item{[JZ1]} W. B. Johnson and M. Zippin, {\sl On subspaces of quotients of
$(\Sigma G)_{\ell _p}$ and $(\Sigma G) _{c_0}$,}\/  {\bf Israel J. Math. 13
nos. 3 and 4} (1972), 311--316.

\item{[JZ2]} W. B. Johnson and M. Zippin,  {\sl Extension of operators from
subspaces of  $c_0(\gamma)$ into $C(K)$ spaces,\/}  {\bf Proc.\ AMS 107 no.
3} (1989), 751--754.

\item{[Lin]} J. Lindenstrauss, {\sl  Extension of compact operators,\/}
{\bf    Memoirs AMS 48} (1964).  

\item{[LP]} J. Lindenstrauss and A. Pe\l czy\'nski, {\sl  Contributions to  
the theory of the classical Banach spaces,\/} {\bf J. Functional Analysis
8}   (1971), 225--249. 

\item{[LR]} J. Lindenstrauss and H. P. Rosenthal, {\sl Automorphisms in
$c_0$, $\ell_1$, and $m$,\/} {\bf Israel J. Math. 7} (1969), 227--239.

\item{[LT1]} J. Lindenstrauss and L. Tzafriri, {\sl Classical Banach spaces
I,   Sequence spaces,\/}  Springer-Verlag, (1977).   

\item{[LT2]} J. Lindenstrauss and L. Tzafriri, {\sl Classical Banach spaces
II,   Function spaces,\/}  Spring\-er-Verlag, (1979).  

\item{[Mac]} G. Mackey, {\sl Note on a theorem of Murray,\/} {\bf Bull. AMS
52} (1046), 322-325.

\item{[Sam1]} D. Samet, {\sl Vector measures are open maps,\/} {\bf
Math. Oper. Res. 9} (1984), 471--474.

\item{[Sam2]} D. Samet, {\sl Continuous selections for vector
measures,\/} {\bf Math. Oper. Res. 12} (1987), 536--543.

\item{[Zip]} M. Zippin, {\sl  A global approach to certain operator
extension problems,\/} Longhorn Notes, {\bf Lecture Notes in Math. 1470}  
Springer-Verlag (1991), 78--84.

\bigskip

\bigskip

\noindent Department of Mathematics, Texas A\&M University, College
Station TX  77843, U.S.A. \qquad  Email address:  johnson@math.tamu.edu

\medskip 

\noindent Institute of Mathematics, The Hebrew University of Jerusalem,
Jerusalem, Israel   
\break Email address: zippin@math.huji.ac.il

\bye